\documentclass[english,12pt,intoc,bibliography=totoc,BCOR10mm,captions=tableheading,fleqn]{article}
\usepackage{lmodern}

\usepackage[T1]{fontenc}
\usepackage[utf8]{inputenc}
\usepackage[a4paper]{geometry}
\usepackage{babel}
\usepackage{prettyref}
\usepackage{units}
\usepackage{mathrsfs}
\usepackage{amsmath}
\usepackage{amsthm}
\usepackage{amssymb}
\usepackage{stackrel}
\usepackage{nomencl}

\providecommand{\makenomenclature}{\makeglossary}
\makenomenclature
\usepackage[unicode=true,
 bookmarks=true,bookmarksnumbered=true,bookmarksopen=true,bookmarksopenlevel=1,
 breaklinks=true,pdfborder={0 0 0},backref=false,colorlinks=false]
 {hyperref}
\hypersetup{
 pdfauthor={Li, Wenwei},
 pdfpagelayout=OneColumn,pdfnewwindow=true,pdfstartview=XYZ,plainpages=false}

\makeatletter
\numberwithin{equation}{section}
\numberwithin{figure}{section}
\newcommand{\lyxaddress}[1]{
\par {\raggedright #1
\vspace{1.4em}
\noindent\par}
}
\theoremstyle{plain}
\newtheorem{thm}{\protect\theoremname}

\usepackage{amssymb}
\usepackage{mathrsfs}
\usepackage{fancyhdr}    
\usepackage{exscale}	
\usepackage{relsize}
\usepackage{mathcomp}  

\usepackage{geometry}

\AtBeginDocument{
  
}

\makeatother

\providecommand{\theoremname}{Theorem}

\begin{document}

\title{\textsf{A Note on the Classification of Permutation Matrix}}

\author{LI Wenwei}

\maketitle

\lyxaddress{\noindent \begin{center}
\emph{liwenwei@ustc.edu}\\
School of Mathematical Science\\
University of Science and Technology of China\\

\par\end{center}}
\begin{abstract}
This paper is concentrated on the classification of permutation matrix
with the permutation similarity relation, mainly about the canonical
form of a permutational similar equivalence class, the cycle matrix
decomposition  and the  factorization of a permutation matrix or
monomial matrix. 

\vspace{0.5cm}

\textbf{Key Words:} permutation matrix, monomial matrix, permutation
similarity, canonical form, cycle matrix decomposition, cycle factorization

\textbf{AMS2000 Subject Classification:} 15A21, 15A36, 
\end{abstract}

$\ $

\section{Introduction}

The incidence matrix of a projective plane of order $n$ is a 0-1
matrix of order $n^{2}+n+1$. Two projective planes will be isomorphic
if the incidence matrix of one projective plane could be transformed
by permuting the rows and/or columns to the incidence matrix of the
other one. After sorting the rows and columns, the incidence matrix
of a projective plane could be reduced in standard form (not unique).
In the reduced form, the incidence matrix could be divided into some
blocks. Most blocks are permutation matrices (refer \cite{Lam1991PP9B}).
If we keep the incidence matrix reduced and keep the position of every
block when permuting the rows and columns of the reduced incidence
matrix, every permutation matrix will be transformed into another
matrix that is permutationally similar to the original one (refer
\cite{liwenwei2014-LR-PP-RT-EATGC} sec. 1.2). 

This paper will focus on the permutational similarity relation and
the classification of the permutation matrices. It will demonstrate
the standard structure of a general permutation matrix, the canonical
form of a permutation similarity class, how to generate the canonical
form and solve some other related issues. The main theorems will be
proved by two different methods (linear algebra method and the combinatorial
method). The number of permutational similarity classes of permutation
matrices of order $n$ will be mentioned on \prettyref{sec:On-the-Number}.
A similar factorization proposition on monomial matrix will be introduced
at the end.

\section{Preliminary\label{sec:Preliminary}}

Let $n$ be a positive integer, $P$ be a square matrix of order $n$.
If $P$ is a binary matrix (every entry in it is either 0 or 1, also
called 0-1 matrix or (0, 1) matrix) and there is a unique ``1” in
its every row and every column,  then $P$ is called a \emph{permutation
matrix }\index{permutation matrix}. If we substitute the units in
a permutation matrix by other non-zero elements, it will be called
a \emph{monomial matrix }\index{monomial matrix} or a \emph{generalized
permutation matrix}. \index{generalized permutation matrix}  

There is a reason for the name ``\emph{permutation matrix}''. A permutation
matrix $P$ of order $n$ multiplies a matrix $T$ of size $n$ $\times$
$r$ (from the left side of $T$) will result the permutation of rows
of $T$. If $U$ is a matrix of size $t$ by $n$, that $P$ acts
on $U$ from the right hand side of $U$ will leads to the permutation
of columns of $U$. \label{fn:-permutation-matrix}

Let $k$ be a positive integer greater than 1, $C$ be an invertible
(0, 1) matrix of order $k$,  if $C^{k}$ = $I_{k}$ ($I_{k}$ is
the identity matrix of order $k$) and $C^{i}$ $\neq$ $I_{k}$ for
any $i$ ($1<i<k$), then $C$ will be called a \emph{cycle matrix
}\index{cycle matrix} of order $k$. A cycle matrix of order $k$
in this form $\left[\begin{array}{ccccc}
0 &  &  &  & 1\\
1 & 0\\
 & 1 & \ddots\\
 &  & \ddots & \ddots\\
 &  &  & 1 & 0
\end{array}\right]$ will be called a \emph{standard cycle matrix}. \index{cycle matrix (standard)}
Sometimes, the identity matrix of order 1 can also be considered as
a cycle matrix of order 1. 

If $C_{1}$ is a permutation matrix of order $n$, and there are exact
$k$ entries in diagonal being 0, (here $2\leqslant k\leqslant n$),
if $C^{k}$ = $I_{n}$ and $C^{i}$ $\neq$ $I_{n}$ for any $i$
($1\leqslant i<k$), then $C_{1}$ will be called a \emph{Generalized
Cycle Matrix of Type }I \index{Generalized Cycle Matrix of Type I@\emph{Generalized Cycle Matrix of Type }I}
with cycle order $k$. \index{cycle order}

If $C_{2}$ is a (0, 1) matrix of order $n$, rank $C_{2}$ = $k$,
and there are exact $k$ entries in $C_{2}$ is non-zero, ($2\leqslant k\leqslant n$),
if $C^{k}$ is a diagonal of rank $k$, too, and $C^{i}$ is non-diagonal
 ($1\leqslant i<k$), then $C_{2}$ will be called a \emph{Generalized
Cycle Matrix of Type }II \index{Generalized Cycle Matrix of Type II@\emph{Generalized Cycle Matrix of Type} II}
with cycle order $k$. Obviously, a generalized cycle matrix of type
II plus some suitable diagonal (0, 1) matrix will result a Generalized
Cycle Matrix of Type I with the same cycle order. 

Let $A$ and $B$ be two monomial matrix of order $n$, if there is
a permutation matrix $T$ such that $B=T^{-1}AT$, then $A$ and $B$
will be called \emph{permutationally similar}. \index{permutationally similar@\emph{permutationally similar}}
Of course the permutation similarity relation is an equivalence relation.
Hence the set of the permutation matrices (or monomial matrix) or
order $n$ will be split into some equivalence classes.

\section{Main Result\label{sec:main-result}}

There are some problems presented naturally:
\begin{enumerate}
\item What's the canonical form of a permutation similarity class?
\item How to generate the canonical form of a given permutation matrix?
\item If $B$ is the canonical form of the permutation matrix $A$, how
to find the permutation matrix $T$, such that $B=T^{-1}AT$?\end{enumerate}
\begin{thm}
\emph{\label{thm:Decomposition} (}Decomposition Theorem\emph{) }For
any permutation matrix $A$ of order $n$, if $A$ is not identical,
then there are some generalized cycle matrices $Q_{1}$, $Q_{2}$,
$\cdots$, $Q_{r}$ of type \emph{II }and a diagonal matrix $D_{t}$
of rank $t$, such that, $A$ \emph{$=$} $Q_{1}+Q_{2}+\cdots+Q_{r}+D_{t}$,
where the non-zero elements in $D_{t}$ are all ones \emph{, }$\sum\limits _{i=1}^{r}\textrm{\emph{rank}}Q_{i}$
\emph{+} $t$ \emph{= $n$; }$1\leqslant r\leqslant\left\lfloor \dfrac{n}{2}\right\rfloor $,
r, $P_{i}$\emph{ }\textup{(}$i$ \emph{=} $1$, $2$, $\cdots$,
$r$\textup{)}\emph{ }and $D_{t}$ are determined by $A$. 
\end{thm}
If the cycle order of $Q_{i}$ is $k_{i}$,\emph{ }($i$ = $1$, $2$,
$\cdots$, $r$)\emph{, }then $2\leqslant\sum\limits _{i=1}^{r}k_{i}\leqslant n$.
When $A$ is a cycle matrix, $t=0$, $r=1$, $k_{1}$ = $n$. 
\begin{thm}
\emph{\label{thm:Factorization} (}Factorization Theorem\emph{) }For
any permutation matrix $A$ of order $n$, if $A$ is not identical,
then there are some generalized cycle matrices $P_{1}$, $P_{2}$,
$\cdots$, $P_{r}$ of type \emph{I }, such that, $A$ = $P_{1}P_{2}\cdots P_{r}$,
where $1\leqslant r\leqslant\left\lfloor \dfrac{n}{2}\right\rfloor $
\emph{; }r, $P_{i}$ \textup{(}$i$ \emph{=} $1$, $2$, $\cdots$,
$r$\textup{)}\emph{ }are determined by $A$. $P_{i_{1}}$ and $P_{i_{2}}$
commute \emph{($1\leqslant i_{1}\neq i_{2}\leqslant r$) }. 
\end{thm}
If the cycle order of $P_{i}$ is $k_{i}$, ($i$ = $1$, $2$, $\cdots$,
$r$)\emph{, }then $2\leqslant\sum\limits _{i=1}^{r}k_{i}\leqslant n$.
 
\begin{thm}
\emph{\label{thm:Similarity} (}Similarity Theorem\emph{) }For any
permutation matrix $A$ of order $n$, there is a permutation matrix
$T$, such that, $T^{-1}AT$ \emph{=} \emph{$\textrm{diag}$ }$\{I_{t}$,
$N_{k_{1}}$, $\cdots$, $N_{k_{r}}\}$, where $N_{k_{i}}$ \emph{=}
$\left[\begin{array}{ccccc}
0 &  &  &  & 1\\
1 & 0\\
 & 1 & \ddots\\
 &  & \ddots & \ddots\\
 &  &  & 1 & 0
\end{array}\right]$ is a cycle matrix of order $k_{i}$ in standard form,\emph{ }\textup{(}$i$
\emph{=} $1$, $2$, $\cdots$, $r$\textup{)}\emph{, $2$ $\leqslant$
$k_{1}$ $\leqslant$ $k_{2}$ $\cdots$ $\leqslant$ }$k_{r}$, $0\leqslant r\leqslant\left\lfloor \dfrac{n}{2}\right\rfloor $,
$0\leqslant t\leqslant n$, and $\sum\limits _{i=1}^{r}k_{i}$ \emph{+}
$t$ \emph{=} $n$. $T$, t, r, $k_{r}$ are determined by $A$. 
\end{thm}
If $A$ is a identity matrix, then $t$ = $n$, $r$ = 0. When $A$
is a cycle matrix, $t=0$, $r=1$, $k_{1}$ = $n$. In this theorem,
the quasi-diagonal matrix (or block-diagonal matrices) \emph{$\textrm{diag}$
}$\{I_{t}$, $N_{k_{1}}$, $\cdots$, $N_{k_{r}}\}$ will be called
the \emph{canonical form} of a permutation matrix in permutational
similarity relation.

\section{Proof\label{sub:Fit_C1(n)}}

Here a proof by linear algebra method is presented.

\subsubsection*{Proof of \prettyref{thm:Similarity} }

For any permutation matrix $A$ of order $n$, let $\mathscr{A}$
be a linear transformation defined on the vector space $\mathbb{R}^{n}$
with bases ${\cal B}$ = \{$e_{1}$, $e_{2}$, $\cdots$, $e_{n}$\},
where $e_{i}$ = $(\underbrace{0\ \cdots\ 0}_{i-1},\ 1,\ \underbrace{0\ \cdots\ 0}_{n-i})^{\mathrm{T}}$,
\footnote{$\ $ Here the regular letter ``T'' in the upper index means transposition. }
($i$ = 1, 2, $\cdots$, $n$), such that $A$ is the matrix of the
transformation $\mathscr{A}$ in the basis ${\cal B}$, or for any
vector $\alpha\in\mathbb{R}^{n}$ with the coordinates $x$ (in the
basis ${\cal B}$), the coordinates of $\mathscr{A}\alpha$ is $Ax$,
i.e., $\mathscr{A}\alpha={\cal B}Ax$. Here the coordinates are written
in column vectors. 

It is clear that the coordinates of $e_{i}$ in the basis ${\cal B}$
is $(\underbrace{0\ \cdots\ 0}_{i-1},\ 1,\ \underbrace{0\ \cdots\ 0}_{n-i})^{\mathrm{T}}$.
Since $A$ is a permutation matrix, $Ae_{i}$ is the $i$'th column
of $A$. 

Let $S$ = \{1, 2, $\cdots$, $n$\}, ${\cal C}$ = $\left\{ e_{i}\,\bigr|\,i\in S\right\} $,
 $a_{11}$ = min $S$, $F_{1}$ = $\left[a_{11}\right]$, $G_{1}=\left[e_{a_{11}}\right]$.
(Here $F_{1}$ and $G_{1}$ are sequences, or sets equipped with orders)
Of course $Ae_{a_{11}}\in{\cal C}$. If $Ae_{a_{11}}\neq e_{a_{11}}$,
assume $e_{a_{12}}=Ae_{a_{11}}$, then put $a_{12}$ and $e_{a_{12}}$
to the ends of the sequences $F_{1}$ and $G_{1}$, respectively.
If $Ae_{a_{1,j}}\neq e_{a_{11}}$, suppose $e_{a_{1,(j+1)}}=Ae_{a_{1,j}}$,
(i.e., $A^{j}e_{a_{11}}=e_{a_{1,(j+1)}}$), then add $a_{1,(j+1)}$
and $e_{a_{1,(j+1)}}$ to the ends of sequences $F_{1}$ and $G_{1}$,
respectively ($j$ = 1, 2, $\cdots$ ). Since $Ae_{i}\in{\cal C}$
($\forall i\in S$), there will be a $k'_{1}$ such that $Ae_{a_{1,k'_{1}}}=e_{a_{11}}$
(otherwise the sequence $e_{a_{11}}$, $Ae_{a_{11}}$, $A^{2}e_{a_{11}}$,
$A^{3}e_{a_{11}}$, $\cdots$ will be infinite). Suppose that $k'_{1}$
is the minimal integer satisfying this condition ($1\leqslant k'_{1}\leqslant n$).
It is clear that $A^{k'_{1}}e_{a_{11}}=e_{a_{11}}$, $A^{k'_{1}}e_{a_{1j}}=e_{a_{1j}}$,
(1 $\leqslant$ $j$ $\leqslant$ $k'_{1}$). It is possible that
$k'_{1}$ = 1 or $n$. So, at last $\left|F_{1}\right|$ = $\left|G_{1}\right|$
= $k'_{1}$. Then remove the elements in $G_{1}$ from ${\cal C}$,
remove the elements in $F_{1}$ from $S$. 

Now, if $S$ $\neq$ Ø, let $a_{21}$ = min $S$, $F_{2}$ = $\left[a_{21}\right]$,
$G_{2}=\left[e_{a_{21}}\right]$. It is clear that $Ae_{a_{21}}\in{\cal C}$.
\footnote{$\ $ Otherwise $Ae_{a_{21}}\in G_{1}$, since all the elements in
$G_{1}$ are removed from ${\cal C}$, then there is a $k_{0}$, s.
t. , $A^{k_{0}}e_{a_{11}}=Ae_{a_{21}}$, of course $k_{0}\neq0$,
so, $A^{k_{0}-1}e_{a_{11}}=e_{a_{21}}$ as $A$ is invertible, which
means that $e_{a_{21}}$ = $A^{k_{0}-1}e_{a_{11}}$ is in the set
$G_{1}$. Contradiction. } If $A^{i-1}e_{a_{21}}\neq e_{a_{21}}$, suppose $A^{i-1}e_{a_{21}}=e_{a_{2i}}$
($i$ = 2, 3, $\cdots$ ), then add $a_{2i}$ and $e_{a_{2i}}$ to
the ends of the sequences $F_{2}$ and $G_{2}$, respectively. There
will be a $k'_{2}$, such that $A^{k'_{2}}e_{a_{21}}=e_{a_{21}}$
(let $k'_{2}$ be the minimal integer satisfying this condition. It
is possible that $k'_{2}$ = 1 or $n-k'_{1}$). Obviously, $A^{k'_{2}}e_{a_{2i}}=e_{a_{2i}}$,
(1 $\leqslant$ $i$ $\leqslant$ $k'_{2}$). Then remove the elements
in $G_{2}$ from ${\cal C}$, remove the elements in $F_{2}$ from
$S$. If $S$ $\neq$Ø, continue this step and constructing $F_{3}$,
$F_{4}$, $\cdots$ and $G_{3}$, $G_{4}$, $\cdots$. It will stop
since $n$ is finite. 

Assume that we have $F_{1}$, $F_{2}$, $\cdots$, $F_{u}$ and $G_{1}$,
$G_{2}$, $\cdots$, $G_{u}$, such that $\stackrel[i=1]{u}{\bigcup}F_{i}$
= \{ 1, 2, $\cdots$, $n$ \}, $\stackrel[i=1]{u}{\bigcup}G_{i}$
= \{ $e_{1}$, $e_{2}$, $\cdots$, $e_{n}$ \}, $F_{i}\cap F_{j}$
= $G_{i}\cap G_{j}$ = Ø, ($1\leqslant i\neq j\leqslant u$). There
is a possibility that $u$ = 1 (when $A$ is a cycle matrix of order
$n$) or $n$ (when $A$ is an identity matrix). Sort $F_{1}$, $F_{2}$,
$\cdots$, $F_{u}$ by  candinality, we will have $\left|F'_{1}\right|$
$\leqslant$ $\left|F'_{2}\right|$ $\leqslant$ $\cdots$ $\leqslant$
$\left|F'_{u}\right|$. Then sort $G_{i}$ correspondingly, i.e.,
$G'_{i}$ = $\left\{ e_{x}\,\left|\,x\in F'_{i}\right.\right\} $
($i$ = 1, 2, $\cdots$, $u$). Suppose $\left|F'_{1}\right|$ = $\left|F'_{2}\right|$
= $\cdots$ = $\left|F'_{t}\right|$ = 1. It is possible that $t$
= $n$ (when $A$ is an identity matrix) or 0. Let $r=u-t$. Denote
the unique element in $G'_{i}$ by $e'_{i}$ ($i$ = 1, 2, $\cdots$,
$t$). Let $k_{j}$ = $\left|G'_{t+j}\right|$, and denote the elements
in $G'_{t+j}$ by $e'_{j,v}$ ($j$ = 1, 2, $\cdots$, $r$; $v$
= 1, 2, $\cdots$, $k_{j}$). Hence, the matrix of the restriction
of $\mathscr{A}$ in the subspace spanned by the bases ${\cal D}_{0}$
= \{ $e'_{1}$, $e'_{2}$, $\cdots$, $e'_{t}$ \} is $I_{t}$ since
$\mathscr{A}e'_{i}=e'_{i}$ ($i$ = 1, 2, $\cdots$, $t$), or $\mathscr{A}{\cal D}_{0}$
= ${\cal D}_{0}I_{t}$; and the matrix of the restriction of $\mathscr{A}$
in the subspace spanned by the bases ${\cal D}_{j}$ = \{ $e'_{j,1}$,
$e'_{j,2}$, $\cdots$, $e'_{j,k_{j}}$ \} ($j$ = 1, 2, $\cdots$,
$r$) is $N_{k_{i}}$ = $\left[\begin{array}{ccccc}
0 &  &  &  & 1\\
1 & 0\\
 & 1 & \ddots\\
 &  & \ddots & \ddots\\
 &  &  & 1 & 0
\end{array}\right]$, a cycle matrix of order $k_{j}$, as $\mathscr{A}e'_{j,v}=e'_{j,v+1}$,
($v$ = 1, 2, $\cdots$, $k_{j}-1$), and $\mathscr{A}e'_{j,k_{j}}=e'_{j,1}$,
i.e., $\mathscr{A}{\cal D}_{j}$ = ${\cal D}_{j}N_{j}$ ($j$ = 1,
2, $\cdots$, $r$). So, the matrix of $\mathscr{A}$ with bases ${\cal D}$
= \{ $e'_{1}$, $e'_{2}$, $\cdots$, $e'_{t}$, $e'_{1,1}$, $e'_{1,2}$,
$\cdots$, $e'_{1,k_{1}}$, $\bullet\bullet\bullet$ , $e'_{r,1}$,
$e'_{r,2}$, $\cdots$, $e'_{r,k_{r}}$\} is $B$ = $I_{t}$ $\oplus$
$N_{k_{1}}$ $\oplus$ $\cdots$ $\oplus$ $N_{k_{r}}$ = \emph{$\mathrm{diag}$
}$\{I_{t}$, $N_{k_{1}}$, $\cdots$, $N_{k_{r}}\}$. 

Since ${\cal D}$ is a reordering of ${\cal B}$, so there is a permutation
matrix $T$, such that ${\cal D}$ = ${\cal B}T$. Then $B=T^{-1}AT$.
 So  \prettyref{thm:Similarity} is proved.

\subsubsection*{Proof of \prettyref{thm:Decomposition}}

With $F'_{i}$ generated above, construct a 0-1 matrix $D_{t}$ of
order $n$, such that the $j$'th column of $D_{t}$ is the $j$'th
column of $A$ ($j\in\stackrel[i=1]{t}{\bigcup}F'_{i}$) and the other
columns of $D_{t}$ are 0 vectors. Of course, $D_{t}$ is a diagonal
matrix of rank $t$, as the $j$'th column of $A$ is $e_{j}$  (by
definition, $Ae_{j}=e_{j}$). Construct a 0-1 matrix $Q_{i}$ \emph{$($
}$i$ = $1$, $2$, $\cdots$, $r$ \emph{$)$} of order $n$, such
that the $j$'th column of $Q_{i}$ is the $j$'th column of $A$
($j\in F'_{t+i}$) and the other columns of $Q_{i}$ are 0 vectors.
As $\left(\stackrel[i=1]{t}{\bigcup}F'_{i}\right)\bigcup$ $\left(\stackrel[i=1]{r}{\bigcup}F'_{t+i}\right)$
= $\stackrel[i=1]{u}{\bigcup}F_{i}$ = \{ 1, 2, $\cdots$, $n$ \},
and $F_{i_{1}}\cap F_{i_{2}}$ = Ø, ($1\leqslant i_{1}\neq i_{2}\leqslant u$),
so every column of $A$ appears exact once in a  matrix in the expression
$\sum\limits _{i=1}^{r}Q_{i}+D_{t}$ in the same position as it appears
in $A$, besides, the columns in the same position in all the other
addend matrices in the expression $\sum\limits _{i=1}^{r}Q_{i}+D_{t}$
are all 0 vectors,  so, 
\[
\sum\limits _{i=1}^{r}Q_{i}+D_{t}=A.
\]

Now we prove that $Q_{i}$ is a generalized cycle matrix of type II
with cycle order $k_{i}$.

Suppose the members in $F'_{t+i}$ ($i$ = 1, 2, $\cdots$, $r$)
are $a'_{i,1}$, $a'_{i,2}$, $\cdots$, $a'_{i,k_{i}}$. Suppose
$F'_{t+i}=F_{s}$ for some $s$ ($1\leqslant s\leqslant u$), \footnote{$\ $ Then we have a relation about the members in $G'_{t+i}$ and
the members in a certain $G_{s}$, i.e., 
\[
e'_{i,v}=e_{a'_{i,v}},\quad(v=1,2,\cdots,k_{i}),
\]
here $e_{a'_{i,v}}\in G_{s}$. } By the definition of $F_{s}$, we know that $Ae_{a'_{i,v}}$ = $e_{a'_{i,v+1}}$,
($v$ = 1, 2, $\cdots$, $k_{i}-1$), $Ae_{a'_{i,k_{i}}}$ = $e_{a'_{i,1}}$,
so $e_{a'_{i,v+1}}$ is the $a'_{i,v}$'th column of $A$. 

As $Q_{i}$ is made of the 0 vector and $k_{i}$ columns of $A$.
The columns of $A$ are linear independent, so the rank of $Q_{i}$
is $k_{i}$. While the $a'_{i,v}$'th column of $Q_{i}$ is the $a'_{i,v}$'th
column of $A$, so $Q_{i}e_{a'_{i,v}}$ = $e_{a'_{i,v+1}}$, ($v$
= 1, 2, $\cdots$, $k_{i}-1$), $Q_{i}e_{a'_{i,k_{i}}}$ = $e_{a'_{i,1}}$,
$Q_{i}e_{l}$ = 0 ($\forall e_{l}\in{\cal D}\left\backslash G'_{t+i}\right.$).
Therefore $Q_{i}^{v}e_{a'_{i,1}}$ = $e_{a'_{i,v+1}}$ ($v$ = 1,
2, $\cdots$, $k_{i}-1$), $Q_{i}^{k_{i}}e_{a'_{i,1}}$ = $e_{a'_{i,1}}$,
(so $Q_{i}^{v}$ is not diagonal as $Q_{i}^{v}e_{a'_{i,1}}$ = $e_{a'_{i,v+1}}$
$\neq$ $e_{a'_{i,1}}$). Then $Q_{i}^{k_{i}}e_{a'_{i,v}}$ = $Q_{i}^{v-1}\left(Q_{i}^{k_{i}-v+1}e_{a'_{i,v}}\right)$
= $Q_{i}^{v-1}\left(e_{a'_{i,1}}\right)$ = $e_{a'_{i,v}}$, ($v$
= 1, 2, $\cdots$, $k_{i}$). Hence $Q_{i}^{k_{i}}$ is a diagonal
matrix of rank $k_{i}$. Therefore $Q_{i}$ is a generalized cycle
matrices of type II with cycle order $k_{i}$.

\subsubsection*{Proof of \prettyref{thm:Factorization}}

Let $D_{t}^{(b)}$ = $I_{n}-D_{t}$. Construct a 0-1 matrix $J_{i}^{a}$\emph{
}($i$ = $1$, $2$, $\cdots$, $r$) of order $n$, such that the
$j$'th column of $J_{i}^{a}$ is the $j$'th column of $I_{n}$ ($j\in F'_{t+i}$)
and the other columns of $J_{i}^{a}$ are 0 vectors. Let $J_{i}^{(b)}$
= $I_{n}-J_{i}^{(a)}$. So, $\sum\limits _{i=1}^{r}J_{i}^{(a)}+D_{t}$
= $I_{n}$. It is obvious that $D_{t}^{(b)}D_{t}$ = $D_{t}D_{t}^{(b)}$
= 0, $J_{i}^{(a)}J_{i}^{(b)}$ = $J_{i}^{(b)}J_{i}^{(a)}$ = 0, $Q_{i}J_{i}^{(b)}$
= $J_{i}^{(b)}Q_{i}$ = 0, and $Q_{i_{1}}J_{i_{2}}^{(a)}$ = $J_{i_{2}}^{(a)}Q_{i_{1}}$
= 0, $Q_{i_{1}}J_{i_{2}}^{(b)}$ = $J_{i_{2}}^{(b)}Q_{i_{1}}$ = $Q_{i_{1}}$
$\neq$ 0 ($1\leqslant i_{1}\neq i_{2}\leqslant r$). Although $J_{i_{1}}^{(a)}J_{i_{2}}^{(a)}$
= $J_{i_{2}}^{(a)}J_{i_{1}}^{(a)}$ = 0, but $J_{i_{1}}^{(a)}J_{i_{2}}^{(b)}$
= $J_{i_{2}}^{(b)}J_{i_{1}}^{(a)}$ = $J_{i_{1}}^{(a)}$ $\neq0$.
It is not difficult to prove that $J_{i_{1}}^{(b)}J_{i_{2}}^{(b)}$
= $J_{i_{2}}^{(b)}J_{i_{1}}^{(b)}$ = $I_{n}-J_{i_{2}}^{(a)}-J_{i_{1}}^{(a)}$.

Let $P_{i}=Q_{i}+J_{i}^{(b)}$, so $\mathrm{rank}$ $P_{i}$ = $n$,
$I_{n}+Q_{i}$ = $P_{i}+J_{i}^{(a)}$. \\
\hphantom{let } $P_{i_{1}}P_{i_{2}}$ = $\left(Q_{i_{1}}+I_{n}-J_{i_{1}}^{(a)}\right)$
$\left(Q_{i_{2}}+I_{n}-J_{i_{2}}^{(a)}\right)$ = $Q_{i_{1}}+Q_{i_{2}}+I_{n}-J_{i_{1}}^{(a)}-J_{i_{2}}^{(a)}$,
then \emph{
\[
\prod\limits _{i=1}^{r}P_{i}=\prod\limits _{i=1}^{r}\left(Q_{i}+I_{n}-J_{i}^{(a)}\right)=\sum\limits _{i=1}^{r}Q_{i}+I_{n}-\sum\limits _{i=1}^{r}J_{i}^{(a)}=\sum\limits _{i=1}^{r}Q_{i}+D_{t}=A.
\]
}

We can prove the equality  above in another way. 

It is clear that $Q_{i_{1}}Q_{i_{2}}=0$, $Q_{i_{1}}D_{t}$ = $D_{t}Q_{i_{1}}$
= 0 ($1\leqslant i_{1}\neq i_{2}\leqslant r$). So 
\begin{equation}
\left(I_{n}+D_{t}\right)\prod\limits _{i=1}^{r}\left(I_{n}+Q_{i}\right)=I_{n}+\sum\limits _{i=1}^{r}Q_{i}+D_{t}=I_{n}+A.\label{eq:Fml-1A}
\end{equation}

Since $D_{t}Q_{i}$ = $Q_{i}D_{t}$ = 0, so $D_{t}J_{i}^{(a)}$ =
$J_{i}^{(a)}D_{t}$ = 0. By definition, when $1\leqslant i\leqslant r$,
the $v$'th column or the $v$'th row of  $\prod\limits _{\begin{array}{c}
_{1\leqslant j\leqslant r}\\
^{j\neq i}
\end{array}}^{r}P_{j}$ with $v\in F'_{t+i}\bigcup\left(\stackrel[i=1]{t}{\bigcup}F'_{i}\right)$
are same as the $v$'th column or the $v$'th row of $I_{n}$, respectively,
(since the $v$'th column or the $v$'th row of $P_{j}$ ($1\leqslant j\leqslant r$,
$j\neq i$) with $v\in F'_{t+i}\bigcup\left(\stackrel[i=1]{t}{\bigcup}F'_{i}\right)$
are the same as the $v$'th column or the $v$'th row of $I_{n}$),
so when $J_{i}^{(a)}$ is multiplied by $\prod\limits _{\begin{array}{c}
_{1\leqslant j\leqslant r}\\
^{j\neq i}
\end{array}}^{r}P_{j}$ (no matter from the left or right), the $v$'th column and the the
$v$'th row will not change, while the other columns and rows of $J_{i}^{(a)}$
is 0 vector, so $J_{i}^{(a)}\prod\limits _{\begin{array}{c}
_{1\leqslant j\leqslant r}\\
^{j\neq i}
\end{array}}^{r}P_{j}$ = $J_{i}^{(a)}$. For the same reason, $D_{t}\prod\limits _{i=1}^{r}P_{i}$
= $D_{t}$. As

\phantom{$=$} $\ $ $\left(I_{n}+D_{t}\right)\prod\limits _{i=1}^{r}\left(I_{n}+Q_{i}\right)$
= $\left(I_{n}+D_{t}\right)\prod\limits _{i=1}^{r}\left(P_{i}+J_{i}^{(a)}\right)$
\\
= $\left(I_{n}+D_{t}\right)\left(\prod\limits _{i=1}^{r}P_{i}+\sum\limits _{i=1}^{r}\left(J_{i}^{(a)}\prod\limits _{\begin{array}{c}
_{1\leqslant j\leqslant r}\\
^{j\neq i}
\end{array}}^{r}P_{j}\right)\right)$ = $\left(I_{n}+D_{t}\right)\left(\prod\limits _{i=1}^{r}P_{i}+\sum\limits _{i=1}^{r}J_{i}^{(a)}\right)$
\\
= $\prod\limits _{i=1}^{r}P_{i}+\sum\limits _{i=1}^{r}J_{i}^{(a)}+D_{t}\prod\limits _{i=1}^{r}P_{i}+D_{t}\sum\limits _{i=1}^{r}J_{i}^{(a)}$
= $\prod\limits _{i=1}^{r}P_{i}+\sum\limits _{i=1}^{r}J_{i}^{(a)}+D_{t}+0$
= $\prod\limits _{i=1}^{r}P_{i}+I_{n}$, 

that is, 
\begin{equation}
\left(I_{n}+D_{t}\right)\prod\limits _{i=1}^{r}\left(I_{n}+Q_{i}\right)=\prod\limits _{i=1}^{r}P_{i}+I_{n}.\label{eq:Fml-1B}
\end{equation}

By equations \eqref{eq:Fml-1A} and \prettyref{eq:Fml-1B}, we will
have $\prod\limits _{i=1}^{r}P_{i}+I_{n}=I_{n}+A$, \\
so $\prod\limits _{i=1}^{r}P_{i}=A$. 

Now we prove that $P_{i}$ is a generalized cycle matrix of type II
with cycle order $k_{i}$. 

Since $P_{i}=Q_{i}+J_{i}^{(b)}$, $Q_{i}J_{i}^{(b)}$ = $J_{i}^{(b)}Q_{i}$
= 0, so $P_{i}^{m}$ = $Q_{i}^{m}+\left(J_{i}^{(b)}\right)^{m}$ =
$Q_{i}^{m}+J_{i}^{(b)}$ ($\forall m\in\mathbb{Z}^{+}$), then $P_{i}^{k_{i}}=Q_{i}^{k_{i}}+J_{i}^{(b)}$.
As $Q_{i}^{k_{i}}e_{a'_{i,v}}$ = $A^{k_{i}}e_{a'_{i,v}}$ = $e_{a'_{i,v}}$,
($v$ = 1, 2, $\cdots$, $k_{i}$); $Q_{i}e_{l}$ = 0 $\Longrightarrow$
$Q_{i}^{k_{i}}e_{l}$ = 0 ($\forall e_{l}\in{\cal D}\left\backslash G'_{t+i}\right.$).
On the other hand, $J_{i}^{(b)}e_{a'_{i,v}}$ = 0 ($v$ = 1, 2, $\cdots$,
$k_{i}$), $J_{i}^{(b)}e_{l}$ = $e_{l}$, ($\forall e_{l}\in{\cal D}\left\backslash G'_{t+i}\right.$).
So for any $e_{l}$ in ${\cal B}$, if $e_{l}\in G'_{i}$, $\left(Q_{i}^{k_{i}}+J_{i}^{(b)}\right)e_{l}$
= $Q_{i}^{k_{i}}e_{l}$ = $e_{l}$, if $e_{l}\notin G'_{i}$, $\left(Q_{i}^{k_{i}}+J_{i}^{(b)}\right)e_{l}$
= $J_{i}^{(b)}e_{l}$ = $e_{l}$. That means $P_{i}^{k_{i}}e_{l}$
= $\left(Q_{i}^{k_{i}}+J_{i}^{(b)}\right)e_{l}$ = $e_{l}$ ($\forall e_{l}\in{\cal B}$).
So \hphantom{so } $P_{i}^{k_{i}}$ $\left(e_{1},e_{2},\cdots,e_{n}\right)$
= $\left(e_{1},e_{2},\cdots,e_{n}\right)$, or $P_{i}^{k_{i}}$ $I_{n}$
= $I_{n}$, hence $P_{i}^{k_{i}}$ = $I_{n}$. (Actually, $Q_{i}^{k_{i}}$
= $J_{i}^{(a)}$, so $P_{i}^{k_{i}}$ = $Q_{i}^{k_{i}}+J_{i}^{(b)}$
= $J_{i}^{(a)}+J_{i}^{(b)}$ = $I_{n}$.) When $1\leqslant m<k_{i}$,
$Q_{i}^{m}$ is not diagonal, neither is $Q_{i}^{m}+J_{i}^{(b)}$
= $P_{i}^{m}$. So $P_{i}$ is a generalized cycle matrix of type
II with cycle order $k_{i}$. 

$\ $

For instance, for the matrix $P_{1}$ = $\left[\begin{array}{cccccc}
0 & 0 & 0 & 0 & 1 & 0\\
0 & 0 & 0 & 1 & 0 & 0\\
0 & 1 & 0 & 0 & 0 & 0\\
0 & 0 & 1 & 0 & 0 & 0\\
1 & 0 & 0 & 0 & 0 & 0\\
0 & 0 & 0 & 0 & 0 & 1
\end{array}\right]$, $P_{1}$ $e_{1}$ = $e_{5}$, $P_{1}$ $e_{5}$ = $e_{1}$, so $F_{1}$
= {[}1, 5{]}, $\left|F_{1}\right|=2$; $P_{1}$ $e_{2}$ = $e_{3}$,
$P_{1}$ $e_{3}$ = $e_{4}$, $P_{1}$ $e_{4}$ = $e_{2}$, so $F_{2}$
= {[}2, 3, 4{]}, $\left|F_{2}\right|=3$; $P_{1}$ $e_{6}$ = $e_{6}$,
$F_{3}$ = {[}6{]}, $\left|F_{3}\right|=1$. So $P_{1}$ is permutationaly
similar to the canonical form $B_{1}$ = \emph{$\textrm{diag}$ }$\{I_{1}$,
$N_{2}$, $N_{3}\}$ = $\left[\begin{array}{cccccc}
1\\
 & 0 & 1\\
 & 1 & 0\\
 &  &  & 0 & 0 & 1\\
 &  &  & 1 & 0 & 0\\
 &  &  & 0 & 1 & 0
\end{array}\right]$, or

$P_{1}$ $\left\{ e_{6};\,e_{1},e_{5};\,e_{2},e_{3},e_{4}\right\} $
= $\left\{ e_{6};\,e_{5},e_{1};\,e_{3},e_{4},e_{2}\right\} $ = $\left\{ e_{6};\,e_{1},e_{5};\,e_{2},e_{3},e_{4}\right\} $
$\left[\begin{array}{cccccc}
1\\
 & 0 & 1\\
 & 1 & 0\\
 &  &  & 0 & 0 & 1\\
 &  &  & 1 & 0 & 0\\
 &  &  & 0 & 1 & 0
\end{array}\right]$.

As $\left\{ e_{6};\,e_{1},e_{5};\,e_{2},e_{3},e_{4}\right\} $ = $\left\{ e_{1},e_{2},e_{3},e_{4},e_{5},e_{6}\right\} $
$\left[\begin{array}{cccccc}
0 & 1 & 0\\
 &  & 0 & 1\\
 &  & 0 &  & 1\\
 &  & 0 &  &  & 1\\
 &  & 1\\
1 & 0 & 0
\end{array}\right]$, let $T_{1}$ = $\left[\begin{array}{cccccc}
0 & 1 & 0\\
 &  & 0 & 1\\
 &  & 0 &  & 1\\
 &  & 0 &  &  & 1\\
 &  & 1\\
1 & 0 & 0
\end{array}\right]$, then $T_{1}^{-1}$ = $T_{1}^{\mathrm{T}}$ = $\left[\begin{array}{cccccc}
0 &  &  &  &  & 1\\
1 &  &  &  &  & 0\\
0 & 0 & 0 & 0 & 1 & 0\\
 & 1\\
 &  & 1\\
 &  &  & 1
\end{array}\right]$, so $\left[\begin{array}{cccccc}
0 & 0 & 0 & 0 & 1 & 0\\
0 & 0 & 0 & 1 & 0 & 0\\
0 & 1 & 0 & 0 & 0 & 0\\
0 & 0 & 1 & 0 & 0 & 0\\
1 & 0 & 0 & 0 & 0 & 0\\
0 & 0 & 0 & 0 & 0 & 1
\end{array}\right]$ = $\left[\begin{array}{cccccc}
0 & 1 & 0\\
 &  & 0 & 1\\
 &  & 0 &  & 1\\
 &  & 0 &  &  & 1\\
 &  & 1\\
1 & 0 & 0
\end{array}\right]$ $\left[\begin{array}{cccccc}
1\\
 & 0 & 1\\
 & 1 & 0\\
 &  &  & 0 & 0 & 1\\
 &  &  & 1 & 0 & 0\\
 &  &  & 0 & 1 & 0
\end{array}\right]$ $\left[\begin{array}{cccccc}
0 &  &  &  &  & 1\\
1 &  &  &  &  & 0\\
0 & 0 & 0 & 0 & 1 & 0\\
 & 1\\
 &  & 1\\
 &  &  & 1
\end{array}\right]$, ( $P_{1}$ = $T_{1}$\emph{ }$B_{1}$ $T_{1}^{-1}$).

$\left[\begin{array}{cccccc}
1\\
 & 0 & 1\\
 & 1 & 0\\
 &  &  & 0 & 0 & 1\\
 &  &  & 1 & 0 & 0\\
 &  &  & 0 & 1 & 0
\end{array}\right]$ = $\left[\begin{array}{cccccc}
0 &  &  &  &  & 1\\
1 &  &  &  &  & 0\\
0 & 0 & 0 & 0 & 1 & 0\\
 & 1\\
 &  & 1\\
 &  &  & 1
\end{array}\right]$ $\left[\begin{array}{cccccc}
0 & 0 & 0 & 0 & 1 & 0\\
0 & 0 & 0 & 1 & 0 & 0\\
0 & 1 & 0 & 0 & 0 & 0\\
0 & 0 & 1 & 0 & 0 & 0\\
1 & 0 & 0 & 0 & 0 & 0\\
0 & 0 & 0 & 0 & 0 & 1
\end{array}\right]$ $\left[\begin{array}{cccccc}
0 & 1 & 0\\
 &  & 0 & 1\\
 &  & 0 &  & 1\\
 &  & 0 &  &  & 1\\
 &  & 1\\
1 & 0 & 0
\end{array}\right]$, ($B_{1}=T_{1}^{-1}P_{1}T_{1}$).

\section{On the Number of Permutation Similarity Classes\label{sec:On-the-Number}}

The number of permutation similarity classes of permutation matrices
of order $n$ is the partition number $p(n)$. There is a recursion
for $p(n)$, 
\begin{align}
p(n) & =p(n-1)+p(n-2)-p(n-5)-p(n-7)+\cdots+\nonumber \\
 & \ \ \ \ (-1)^{k-1}p\left(n-\dfrac{3k^{2}\pm k}{2}\right)+\cdots\cdots\nonumber \\
 & =\sum\limits _{k=1}^{k_{1}}(-1)^{k-1}p\left(n-\dfrac{3k^{2}+k}{2}\right)+\sum\limits _{k=1}^{k_{2}}(-1)^{k-1}p\left(n-\dfrac{3k^{2}-k}{2}\right),\label{eq:pn-recursion}
\end{align}
(Refer \cite{Marshall1958Survey}, page 55), where 
\begin{equation}
k_{1}=\left\lfloor \dfrac{\sqrt{24n+1}-1}{6}\right\rfloor ,\ k_{2}=\left\lfloor \dfrac{\sqrt{24n+1}+1}{6}\right\rfloor ,\label{eq:pn-recursion-k1-k2}
\end{equation}
 and assume that $p(0)=1$. Here $\left\lfloor x\right\rfloor $ is
the floor function, it stands for the maximum integer that is less
than or equal to the real number $x$.

We may find a famous asymptotic formula for $p(n)$ in references
\cite{Eric1999PtFuncP} or \cite{NIST2015FNTANTUP}, \label{Sym:Rh(n)}
\nomenclature[Rh(n)]{$R_{\mathrm{h}}(n)$}{The Hardy-Ramanujan's asymptotic formula. \pageref{Sym:Rh(n)}}
\begin{equation}
p(n)\sim\dfrac{1}{4n\sqrt{3}}\exp\left(\sqrt{\frac{2}{3}}\pi n^{\nicefrac{1}{2}}\right).\label{eq:Ram-Hardy-Fml-Est-1}
\end{equation}
 This formula is obtained by Godfrey H. Hardy and Srinivasa Ramanujan
in 1918 in the famous paper \cite{Ramanujan1918AsymFmlComAnal}. (In
\cite{Pal1942ElemProofRHFml} and \cite{Newman1962SimProofPtFml},
we can find two different proofs of this formula. The evaluation of
the constants can be found in \cite{Donald1951EvalfConstHRFml}.)

Formula \ref{eq:Ram-Hardy-Fml-Est-1} is very import for analysis
in theory. It is very convenient to estimate the value of $p(n)$
especially for ordinary people not majored in mathematics. \footnote{$\ $ Compared with another famous formula in convergent series found
by Rademacher in 1937, based on the work of Hardy and Srinivasa Ramanujan,
refer \cite{Marshall1958Survey} or \cite{Rademacher1937ConvergSeries}. } But the accuracy is not so satisfying when $n$ is small. 

In \cite{liwenwei2016-Estmn-pn-arXiv}, several other formulae modified
from formula \prettyref{eq:Ram-Hardy-Fml-Est-1} is obtained (with
high accuracy). Such as 
\begin{equation}
p(n)\approx\left\lfloor \dfrac{\exp\left(\sqrt{\frac{2}{3}}\pi\sqrt{n}\right)}{4\sqrt{3}\left(n+C'_{2}(n)\right)}+\dfrac{1}{2}\right\rfloor ,\quad1\leqslant n\leqslant80.
\end{equation}
 with a relative error less than 0.004\%, where 
\[
C'_{2}(n)=\begin{cases}
0.4527092482\times\sqrt{n+4.35278}-0.05498719946, & n=3,5,7,\cdots,79;\\
0.4412187317\times\sqrt{n-2.01699}+0.2102618735, & n=4,6,8\cdots,80.
\end{cases}
\]
 and 
\begin{equation}
p(n)\approx\left\lfloor \dfrac{\exp\left(\sqrt{\frac{2}{3}}\pi\sqrt{n}\right)}{4\sqrt{3}\left(n+a_{2}\sqrt{n+c_{2}}+b_{2}\right)}+\dfrac{1}{2}\right\rfloor ,\ \ n\geqslant80
\end{equation}
with a relative error less than $5\times10^{-8}$ when $n\geqslant180$.
Here $a_{2}=0.4432884566$, $b_{2}=0.1325096085$ and $c_{2}=0.274078$.

\section{Result on Monomial Matrix}

For any monomial matrix $M$, it can be written as the product of
a permutation matrix $P$ and an invertible diagonal matrix $D$.
\footnote{$\ $ Turn all the non-zero elements in $M$ into 1, then we will
have a permutation matrix $P$.  Suppose the unique non-zero elements
in the $i$'th row of $M$ is $c_{i}$, the unique non-zero elements
in the $i$'th column of $M$ is $d_{i}$, $i$ = 1, 2, $\cdots$,
$n$. Let $D_{1}$ = $\mathrm{diag}$ \{ $c_{1}$, $c_{2}$, $\cdots$,
$c_{n}$ \}, $D_{2}$ = $\mathrm{diag}$ \{ $d_{1}$, $d_{2}$, $\cdots$,
$d_{n}$ \},  $M$ = $PD_{2}$ = $D_{1}P$. } For the permutation matrix $P$, there is a permutation matrix $T$
such that $T^{-1}PT$ = $Y$ is in canonical form diag\emph{ }$\{I_{t}$,
$N_{1}$, $\cdots$, $N_{r}\}$ as mentioned in \prettyref{thm:Similarity}.
In the expression $T^{-1}PT$, the permutation matrix $T^{-1}$ changes
only the position of the rows, $T$ just changes the position of the
columns, neither will change the values of the members, as the non-zero
members in $M$ and $P$ share the same positions, so do $T^{-1}MT$
and $T^{-1}PT$. Suppose the unique non-zero element in the $i$'th
row of $T^{-1}MT$ is $a_{i}$, the unique non-zero element in the
$i$'th column of $T^{-1}MT$ is $b_{i}$, $i$ = 1, 2, $\cdots$,
$n$. Let $D_{3}$ = $\mathrm{diag}$ \{ $a_{1}$, $a_{2}$, $\cdots$,
$a_{n}$ \}, $D_{4}$ = $\mathrm{diag}$ \{ $b_{1}$, $b_{2}$, $\cdots$,
$b_{n}$ \}, then $T^{-1}MT$ = $D_{3}Y$ = $YD_{4}$.

So $M$ = $D_{1}T\left[\begin{array}{cccc}
I_{t}\\
 & N_{1}\\
 &  & \ddots\\
 &  &  & N_{r}
\end{array}\right]T^{-1}$ = $T\left[\begin{array}{cccc}
I_{t}\\
 & N_{1}\\
 &  & \ddots\\
 &  &  & N_{r}
\end{array}\right]T^{-1}D_{2}$\\
\hphantom{So $M$ }= $TD_{3}\left[\begin{array}{cccc}
I_{t}\\
 & N_{1}\\
 &  & \ddots\\
 &  &  & N_{r}
\end{array}\right]T^{-1}$ = $T\left[\begin{array}{cccc}
I_{t}\\
 & N_{1}\\
 &  & \ddots\\
 &  &  & N_{r}
\end{array}\right]D_{4}T^{-1}$.

\section*{Acknowledgements}

\addcontentsline{toc}{section}{Acknowledgements}

The author would like to express the gratitude to his supervisor Prof.\emph{
LI} \emph{Shangzhi }from BUAA (Beihang University in China)  for
his valuable advice. 

\bibliographystyle{amsplain}
\phantomsection\addcontentsline{toc}{section}{\refname}\bibliography{Ref-org-ch-15}

\end{document}